\documentclass[a4paper, 11pt, oneside]{article}

\usepackage[affil-it]{authblk}
\usepackage{blindtext}
\usepackage{amsfonts,amsmath,amssymb,amsxtra,amsthm}
\usepackage{mathrsfs}
\usepackage{times}
\usepackage{anysize}
\usepackage{graphicx}

\input xy
\xyoption{all}
\usepackage[all,knot]{xy}

\newtheorem{theorem}{Theorem}
\newtheorem{corollary}[theorem]{Corollary}

\theoremstyle{definition}

\usepackage{lipsum}

\newcommand{\footnotewithoutnum}[1]{%
  \begingroup\def\thefootnote{}\footnotetext{#1}\endgroup}

\begin{document}

\title{Subword complexes and 2-truncated cubes}

\author{Mikhail Gorsky}
\date{}

\maketitle

Let $W$ be a finite Coxeter group, $S = \left\{s_1,\ldots,s_n\right\}$ be a set of simple reflections generating $W.$ Consider a word $\mathbf{Q}:=\mathbf{q}_1\ldots\mathbf{q}_m$ in the alphabet of simple reflections ($\mathbf{q}_i \in S$ \ $\forall\,i=1,\ldots,m$) and an element $\pi$ of the group $W.$ {\it The subword complex $\Delta(\mathbf{Q};\pi)$} is a pure simplicial complex on the set of vertices $\{\mathbf{q}_1,\ldots,\mathbf{q}_m\}$ corresponding to the letters (more precisely, to their positions) in the word $Q.$ A set of vertices yields a simplex if the complement in $Q$ to the corresponding subword contains a reduced expression of $\pi.$ The maximal simplices correspond to the complements of reduced expressions of $\pi$ in the word $\mathbf{Q}.$ Subword complexes were introduced by A.~Knutson and E.~Miller in the article \cite{1}. They showed in \cite{2} that $\Delta(\mathbf{Q};\pi)$ is spherical if and only if the Demazure product of the word $Q$ equals $\pi;$ otherwise, $\Delta(\mathbf{Q};\pi)$ is a triangulated ball. For spherical subword complexes, there arise natural questions of the existence, of the combinatorial description and of geometric realizations of their polar dual polytopes.

\footnotewithoutnum{By ``polytopes'' we always mean {\bf combinatorial polytopes}. It is known (see \cite{8} and references therein) that each combinatorial 2-truncated cube can be realized as a geometric polytope. By unfortunate abuse of terminology, we speak about \emph{brick polytopes} which were defined in \cite{3} as geometric objects, but we actually mean their combinatorial types.}

In the group $W$ there exists the unique longest element denoted by $w_o.$ We will consider subword complexes of the form $\Delta(\mathbf{c}\mathbf{w}_{o};w_o)$, where $\mathbf{c}$ is a reduced expression of a Coxeter element, $\mathbf{w}_o$ is an arbitrary reduced expression of $w_o$. Such complexes admit a realization by {\it brick polytopes} of V.Pilaud--\allowbreak C.~Stump \cite{3} that we will denote by $\mathbf{B}(\mathbf{c}\mathbf{w}_o;w_o)$. For each Coxeter element $c$ in the group $W$ there are defined the \textit{$c$-cluster complex of type $W$} and its dual polytope -- \textit{(the generalized) $c-$associahedron of type $W$}. These objects are important in the theory of cluster algebras. C.~Ceballos, J.-P.~Labb\'{e} and C.~Stump \cite{4} proved that the complexes $\Delta(\mathbf{c}\mathbf{w}_o(\mathbf{c});w_o)$, where $\mathbf{w}_o(\mathbf{c})$ is the so-called \textit{$\mathbf{c}$-sorting word} for $w_o,$ are the $c$-cluster complexes of type $W.$ Therefore, the polytopes $\mathbf{B}(\mathbf{c}\mathbf{w}_o(\mathbf{c});w_o)$ realize the $c-$associahedra of type $W.$ The choice of a Coxeter element $c$ is equivalent to the choice of a quiver $Q$ being an orientation of the Coxeter diagram of the group $W.$ Let $\mathbf{c}'$ be a reduced expression of an arbitrary (possibly coinciding with $c$) Coxeter element $c'$ in the group $W,$ and let $Q'$ be the corresponding quiver. Let $\widetilde{Q}_{c,c'}$ be the quiver obtained from $Q$ by erasing all edges oriented differently in $Q$ and $Q'$ (we do not remove any vertices); $\widetilde{W}_{c,c'}$ be the Coxeter group, such that $\widetilde{Q}_{c,c'}$ is an orientation of its Coxeter diagram. We can consider $c$ as a Coxeter element in $\widetilde{W}_{c,c'}.$

\begin{theorem}
\label{th1} $\Delta(\mathbf{c}\mathbf{w}_o(\mathbf{c}');w_o)$ is the $c$-cluster complex of type ~$\widetilde{W}_{c,c'}$, ~$\mathbf{B}(\mathbf{c}\mathbf{w}_o(\mathbf{c}');w_o)$ realizes the corresponding generalized associahedron, In particular, let ~$\mathbf{c}_{\rm rev}$ denote the word ~$\mathbf{c}$ written in the opposite direction\footnote{It is a reduced expression of $c^{-1}$.}, then $\mathbf{B}(\mathbf{c}\mathbf{w}_o(\mathbf{c}_{\rm rev});w_o)$ is a (combinatorial) cube.
\end{theorem}

The proof is based on arguments similar to those from the article \cite{4}. We use \cite[lemma~2.3]{4} and check that all the results in \cite[Section~5]{4}, except for lemma~5.5 and proposition~5.6, hold for an arbitrary reduced expression~$\mathbf{w}_{o}$, and not only for the $\mathbf{c}$-sorting word.

Consider now an arbitrary reduced expression~$\mathbf{w}_{o}$ of the element~$w_o$. There exists a bijection $\operatorname{Lr}_{c,\mathbf{w}_o}$ between the set of letters in the word 
$\mathbf{c}\mathbf{w}_{o}=c_1 c_2\ldots c_n w_1 w_2\ldots w_N$
and the set $\Phi_{\geqslant -1}=-\Pi \sqcup \Phi_{+}$ of almost positive (i.\,e.\ simple negative and all positive) roots in the root system~$\Phi$, associated to~$W$:
$$
\operatorname{Lr}(c_i)=-\alpha_{c_i}; \qquad
\operatorname{Lr}(w_j)=w_1 w_2\ldots w_{j-1}(\alpha_{w_j}),
$$
where $\alpha_w \in \Phi$ is the root corresponding to $w \in W$.
The choice~$\mathbf{w}_{o}$ yields a total order~$<_{\mathbf{w}_o}$
on~$\Phi_+$: $\alpha <_{\mathbf{w}_o} \beta$, if 
$\operatorname{Lr}^{-1}(\alpha)$ goes in~$\mathbf{w}_{o}$
before~$\operatorname{Lr}^{-1}(\beta)$. It is known that such an order satisfies the following condition: for any root subsystem of rank~2 in~$\Phi$ with the canonical generators~$\alpha$,~$\beta$ one either has $a_1 \alpha+b_1 \beta
<_{\mathbf{w}_o} a_2 \alpha+b_2 \beta$ provided $a_1 <
a_2$, or $a_1 \alpha+b_1 \beta <_{\mathbf{w}_o} a_2 \alpha+b_2
\beta$ provided $a_1 > a_2$. In other words, the positive roots of any rank 2 subsystem are ordered in one of two natural ways. The choice of ~$\mathbf{c}$ 
also provides a total order on~$\Phi_+$: $<_{\mathbf{c}}\,=\,<_{\mathbf{w}_o(\mathbf{c})}$.
We will say that a root $\gamma \in \Phi_+$ is
\textit{$(c,\mathbf{w}_o)$-stable}, if for any non-commutative rank 2 subsystem $\langle\alpha,\beta\rangle$ containing~$\gamma$ and such that $\gamma \ne \alpha, \beta$, we have
$$
\alpha <_{\mathbf{c}}\beta \Leftrightarrow \alpha<_{\mathbf{w}_o}\beta.
$$
This condition depends on~$c$, but not on~$\mathbf{c}$.
Let  $\operatorname{Stab}(\mathbf{c},\mathbf{w}_o)$ be the set of the
$(c,\mathbf{w}_o)$-stable roots. The main result of this article is the following theorem. 

\begin{theorem}
\label{th2} {\rm(i)} The vertices of $\Delta(\mathbf{c}\mathbf{w}_o;w_o)$
and, equivalently, the facets of
$\mathbf{B}(\mathbf{c}\mathbf{w}_o;w_o)$ are in a one-to-one correspondence with the simple negative and the \textit{$(\mathbf{c},\mathbf{w}_o)$-stable} positive roots in the system~$\Phi$.

{\rm(ii)} Let expressions~$\mathbf{w}_o$,~$\mathbf{w}_o'$ be such that $\operatorname{Stab}(\mathbf{c},\mathbf{w}_o) \subset
\operatorname{Stab}(\mathbf{c},\mathbf{w}_o')$. Then the complex
$\Delta(\mathbf{c}\mathbf{w}_o';w_o)$ can be obtained from the complex $\Delta(\mathbf{c}\mathbf{w}_o;w_o)$
by a sequence of edge subdivisions. Similarly, the polytope $\mathbf{B}(\mathbf{c}\mathbf{w}_o';w_o)$ can be obtained from the polytope $\mathbf{B}(\mathbf{c}\mathbf{w}_o;w_o)$ by a sequence of truncations of faces of codimension ~$2$.\footnote{As combinatorial polytopes.}

\end{theorem}

The proof is based on results of the work~\cite{5}, where it was shown how braid moves in the group $W$ induce compositions of edge subdivisions and the inverse operations on subword complexes. Then we show that there is a pair of roots ~$\alpha$,
$\beta\in\operatorname{Stab}(\mathbf{c},\mathbf{w}_o)$,
$\alpha<_{\mathbf{w}_o}\beta$, $\beta<_{\mathbf{w}_o'}\alpha$,
such there are no other roots from $\operatorname{Stab}(\mathbf{c},\mathbf{w}_o)$ between them (with respect to the order $<_{\mathbf{w}_o}$). It remains to show that, by a sequence of braid moves, one can reorder letters corresponding to the roots in the interval 
$[\alpha,\beta]_{\mathbf{w}_o}$ in such a way that~$\beta$ would go before~$\alpha$.
One can take as a permutation the $\mathbf{c}-$sorting word, then the result follows from the word property of the group $W$.

The class of \textit{$2$-truncated cubes}, i.e. polytopes which can be obtained from the cube of a fixed dimension by a sequence of truncations of faces of codimension $2$, has interesting properties and includes important families of polytopes (cf.~\cite{8}). Each $2$-truncated cube is a \textit{flag} polytope (i.e. any set of its pairwise intersecting facets has a nonempty intersection). By theorems ~\ref{th1} and~\ref{th2}, we get the following statement.

\begin{corollary}
\label{cor1}
Each polytope of the form $\mathbf{B}(\mathbf{c}\mathbf{w}_o;w_o)$
is a $2$-truncated cube and, therefeore, it is flag. In particular, any generalized associahedron is a $2$-truncated cube.
\end{corollary}

\bigskip
{\bf{\large{Acknowledgements.}}} 
This short announcement note was published in Russian Math. Surveys 69 (2014) no. 3. (in the section \emph{Communications of the Moscow Mathematical Society}). In accordance with the rules of the journal, the present version differs from the author's original only in footnotes, acknowledgements and the current address.

This is a summary of Chapter 4.4 of my Ph.D. thesis \cite{6} defended in 2014 at the Steklov Mathematical Institute. The thesis was written in Russian and, by various reasons, its Chapter 4 has not been translated to English, except for the present note; moreover, the latter was not submitted to arXiv earlier. A translation of \cite[Chapter 4]{6}, with detailed proofs, corrections, clarifications, and some new results will appear in the forthcoming work \cite{7}.  

I am indebted to my supervisor Prof. Victor M. Buchstaber for the inspiration in the work and for his support and patience.
I am very grateful to Drew Armstrong, Fr\'ed\'eric Chapoton, Jean-Philippe Labb\'e, Vincent Pilaud, Evgeny Smirnov, Salvatore Stella, Christian Stump, Hugh Thomas, and Vadim Volodin for valuable discussions on subword complexes and on results of this note and of \cite{7}. This work was supported by R\'eseau de Recherche Doctoral en Math\'ematiques de l’\^{I}le
de France and by RFBR (project 14-01-92612).

\medskip

\textsc{Hausdorff Research Institute for Mathematics, Poppelsdorfer Allee 45, 53115 Bonn, Germany.}\par\nopagebreak
\textit{E-mail address}: \texttt{mikhail.gorsky@iaz.uni-stuttgart.de}

\end{document}